\documentclass[11pt]{article}
\usepackage{amsmath}
\usepackage{amsfonts}
\usepackage{amsthm}
\usepackage{verbatim}
\def\N{\mathbb{N}}
\def\F{\mathbb{F}}

\newtheorem{thm}{Theorem}

\newtheorem{coro}{Corollary}
\begin{document}
\centerline{\bf{ On hyperquadratic continued fractions }}
\centerline{\bf{in power series fields over a finite field }}
\vskip 0.5 cm
\centerline{\bf{by}}
\centerline{\bf{A. Lasjaunias}}
\vskip 0.5 cm
(This note is a brief written account of a talk given at the conference on finite fields which was held in Magdeburg in July 2013.)
 \par The first part of this note is a short introduction on continued fraction expansions for certain algebraic power series. 
In the last part, as an illustration, we present a family of algebraic continued fractions of degree 4, including a toy example 
considered about thirty years ago in a pioneer work in this area. 
\par For a general account concerning the background of the matter presented in this note and for more references, 
the reader may consult W. Schmidt's article [8]. 
 
\section{Formal power series over a finite field }

\par The fields of power series over a finite field (or formal numbers) are known to be analogues of the field of real numbers.
 We have the following analogy, where $q$ is a power of the characteristic $p$ of the finite field $\F_q$, having $q$ elements 
and $T$ is a formal indeterminate :

\begin{center}
\begin{tabular}{l c c c c c}
$\pm 1$ &$\quad$ & & $\quad$ & $\mathbb{F}_q^*$\\
$\cap $      &$\quad $ & $\quad$ & $\quad$ & $\cap$  \\
$\mathbb{Z}$ & $\quad$ & $\longleftrightarrow $ & $\quad$ & $\mathbb{F}_q[T]$\\
$\cap $      &$\quad $ & $\quad$ & $\quad$ &  $\cap $\\
$\mathbb{Q}$ & $\quad$ &$\longleftrightarrow $ & $\quad$ & $\mathbb{F}_q(T)$ \\
$\cap $      &$\quad $ & $\quad$ & $\quad$ & $\cap$  \\
$\mathbb{R}$ &$\quad$  &$\quad $ & $\quad$ & $\mathbb{F}_q((1/T))$\\
\end{tabular}
\end{center}
Hence a real number expanded in base $b$  $$x=\sum_{n\leq k }a_nb^n \in \mathbb{R}$$ 
is replaced by a formal power series in $1/T$, with coefficients in a finite field $\F_q$,  
 $$\alpha=\sum_{n\leq k} a_nT^n \in \mathbb{F}_q ((1/T)).$$ 
 In this note the field $\mathbb{F}_q((1/T))$ will be briefly denoted by $\mathbb{F}(q)$. Following this analogy, the field $\mathbb{F}(q)$ is the completion of the field  $\mathbb{F}_q(T)$ for the absolute value $\vert \alpha \vert =\vert T \vert^k$ where $\vert T \vert >1$ is a fixed real number.
\newline As a trivial expansion for a rational element, we have in $\mathbb{F}_p(T)$
$$\qquad 1/(T-1)=T^{-1}+T^{-2}+\dots+T^{-n}+\dots $$
\newline Even though one can consider power series over an arbitrary base field, the case of a finite base field, which we consider here, is of particular importance. Indeed this finiteness implies the following analogue of a classical result.
$$\alpha \in \mathbb{F}_q(T)\iff (a_n)_{n\leq k}\quad \text{ is
eventually periodic.}$$

\section{Continued fractions in $\mathbb{F}(q)$ }

\par Continued fractions in function fields have long been considered; however Baum and Sweet's article [1], in the frame 
of $\F(2)$, is fundamental in the developement of the subject discussed here. 
\par Every $\alpha \in \mathbb{F}(q)$ can be expanded as a continued
  fraction : 
$$\alpha=a_0+1/(a_1+1/(a_2+1/(\dots=[a_0,a_1,a_2,\dots]$$
where the partial quotients $a_i\in \mathbb{F}_q[T]$ and $\deg(a_i)>0$ for $i>0$. As usual 
$\alpha_{n}=[a_n,a_{n+1},\dots]\in \mathbb{F}(q)$ denotes the tail of the expansion.
\newline This expansion is finite if and only if $\alpha \in \mathbb{F}_q(T)$. 
\newline For instance, in $\mathbb{F}_{13}(T)$, we have: 
$$\qquad
\frac{(T^2-1)^4}{2T^7+2T^5+T^3-T}=[7T,10T,5T,-T,9T,11T,T,5T].$$
Moreover, as in the real case, quadratic power series over a finite field have a particular expansion. Indeed we have: $[\mathbb{F}_q(T,\alpha):\mathbb{F}_q(T)]=2 \iff (a_n)_{n\geq 0}$ is eventually periodic. 
\newline We illustrate this with two examples. The first one is the analogue of the golden mean. 
\newline 1) In $\F(p)\quad :$ $\qquad \omega=[T,T,T,\dots,T,\dots]$ implies $\omega=T+1/\omega.$
\newline  
 2) In $\F(11)\quad :$ $\qquad \alpha=[T,2T,3T,T,2T,3T,\dots,T,2T,3T,\dots]$ satisfies the equation 
 $$(6T^2+1)\alpha^2+(5T^3+9T)\alpha+9T^2+10=0.$$
\par In the sequel we will use polynomials directly connected to continued fractions. Given an infinite sequence $(x_n)_{n\geq 1}$ of variables, we
define recursively the sequence of multivariate polynomials 
 $$<\emptyset>=1, \qquad <x_1>=x_1 \quad \text{and}$$ 
$$<x_1,x_2,\dots,x_n>=x_1<x_2,\dots,x_n>+<x_3,\dots,x_n> \quad \text{for} \quad n\geq 2.$$
These polynomials, which are called
continuants, play a fundamental role in the continued fraction
algorithm. Indeed, for a finite continued
fraction, we have the following basic formula :
$$[a_1,a_2,\dots,a_n]=\frac{<a_1,\dots,a_n>}{<a_2,\dots,a_n>}.$$
\par Finally we present an important statement concerning the growth of the partial quotients for algebraic continued fractions (see for instance [8]). We have
$$t(\alpha)=\limsup_{n\geq 0}[\deg(a_{n+1})/\sum_{0\leq k\leq n}\deg(a_k)] \leq d-2,$$
if $\alpha=[a_0,a_1,a_2,\dots]$ is algebraic over $\mathbb{F}_q(T)$ of degree $d\geq 2$.

\section{Hyperquadratic formal power series }

\par In function fields over a finite field, the Frobenius isomorphism plays an important role. The notion of hyperquadratic power series has emerged from the fundamental articles [1] and [7]. Note that such algebraic power series were considered also and independently by Osgood (concurrently with [1]) and Voloch (concurrently with [7]). Both authors were interested in rational approximation, this matter will not be discussed here, but the reader can be referred to [8].  
\par In the following definition $p$ is the characteristic of $\mathbb{F}_q$ and $r=p^t$ where $t\geq 0$ is an integer.
\par We write $\alpha \in
\mathbb{H}_r(q)$, if $\alpha \in \mathbb{F}(q)\backslash
\mathbb{F}_q(T)$ and there exists $(A,B,C,D)\in (\mathbb{F}_q[T])^4$ with 
  $$\alpha=(A\alpha^r+B)/(C\alpha^r+D).$$
Hence $\mathbb{H}_r(q)$ is a subset of algebraic power series. One can show that an algebraic power series of large degree $d$ is rarely hyperquadratic, however we have the following statements:
\newline 1) $\alpha \in \mathbb{H}_1(q) \iff \alpha$ is quadratic over
  $\mathbb{F}_q(T)$. 
\newline 2) $[\mathbb{F}_q(T,\alpha):\mathbb{F}_q(T)]=3 \Rightarrow \alpha \in \mathbb{H}_p(q)$.
\par Unlike quadratic power series, the continued fraction expansion for an hyperquadratic element is generally difficult to 
establish. However, because of a certain proximity with quadratic elements, in some particular cases, these hyperquadratic 
expansions can be explicitly described. We give here three examples :
\newline 
1) In $\mathbb{H}_r(p)\quad :$ $\quad \alpha =[T,T^r,T^{r^2},\dots,T^{r^k},\dots]$ implies $\alpha=T+1/\alpha^r$. 
\newline  
2) In $\mathbb{H}_2(2)\quad :$  $\quad \alpha =[T^2+1,T^3,T,T^5,T^{[3]},\dots,T^{2^n+1},T^{[2^n-1]},\dots]$ implies $\alpha^3+(T^2+1)\alpha^2+T=0.$
\newline 3) In $\mathbb{H}_4(4)\quad :$  $\quad \alpha =[uT,T^{[3]},uT,T^{[15]},\dots,uT,T^{[4^n-1]},\dots]$ implies
$$T^3\alpha^5+(uT^4+T^2+1)\alpha^4+1=0.$$
The first example is trivial considering the action of the Frobenius isomorphism. Note that in this case the growth of the partial
 quotients is easy to check and, according to the inequality stated at the end of the previous section, we obtain $r-1\leq d-2$ and
 this implies $d=r+1$. 
\newline For the second and third example, we use the notation $T^{[m]}=T,T,\dots,T$ ($m$ times). Note that for the second example 
we have $t(\alpha)=1/2$. In the last example 
the field $\F_4$ is $\F_2(u)$ where $u$ satisfies $u^2+u+1=0$. It is also more tricky to prove that this last expansion satisfies the
 given equation (see [3, p. 100]). 

\section{Hyperquadratic power series of degree 4}
\par Since all algebraic power series of degree less than 4 are hyperquadratic, it is natural to try to characterize, among the power series
of degree 4, those which are hyperquadratic. In this direction, in a joint work with A. Bluher [2], we could prove the following.
\begin{thm} Let $p>3$ be a prime number. Let $(A,B,C)\in \F_q(T)^3$ be such that $12A+C^2=0$. Set 
$$P(X)=AX^4+BX^3+CX^2+1\in \F_q(T)[X].$$ Then there exists a
  non-trivial polynomial $H$, such that $P$ divides $H$, of the form
  $$H(X)=UX^{r+1}+VX^r+WX+Z \in \F_q(T)[X],$$
where $r=p\quad \text{ if} \quad p\equiv 1 \mod 3$ and $r=p^2 \quad
\text{ if}\quad p\equiv 2 \mod 3.$
\end{thm}
\par In their article [7], Mills and Robbins introduced a particular algebraic power series of degree 4 in $\F(13)$, for which the continued fraction
expansion had an apparent regular pattern. They tried to conjecture the sequence of partial quotients. In several articles (see [4], [5]
and [6]), starting from this quartic equation and observing that the solution is hyperquadratic, we have developped a method allowing 
us to describe a very wide family of such hyperquadratic expansions, in all odd characteristic and arbitrary agebraic degree,
 including Mills and Robbins example. In the present note, we are interested in a generalization of this quartic element, in all power series
fields $\F(p)$ where $p>3$. The above theorem implies clearly the following.
\begin{coro} Let $p>3$ be a prime. There exists a unique $\alpha \in \F(p)$ solution of
$$P(X)=(9/32)X^4-TX^3+X^2-8/27=0.\quad (RQE)$$
We have $\alpha =(32/9)T-T^{-1}+....\qquad $ and $\quad \alpha \in \mathbb{H}_{p^i}(p)\quad \text{ if}\quad p\equiv i \mod 3$. 
\end{coro}
The original equation, considered by Mills and Robbins, appeared in a much simpler form and this explains its casual introduction.
Indeed, for $p=13$, if $\alpha$ is the root of $(RQE)$ and $\beta$ is defined by $1/\beta(T)=v\alpha(vT)$, where $v^2=5$, 
then $\beta$ is the unique root in $\F(13)$ of the equation: $X^4+X^2-TX+1=0$ (see [7, p. 404]). 
\par We are now going to describe the continued fraction expansion for the root of $(RQE)$ in $\F(p)$, assuming that $p$ is prime
such that $p\equiv 1 \mod 3$. To do so we need to work in a larger frame. 
\par Here $p$ is an odd prime and $k$ is an integer, with $1\leq k<p/2$. We define in $(\F_p^*)^{2k}$
    the $2k$-tuple $(v_1,v_2,\dots,v_{2k})$ such that $v_1=2k-1$ and
$$v_{i+1}v_i=\frac{(2k-2i-1)(2k-2i+1)}{i(2k-i)} \qquad \text{ for
}\quad 1\leq i\leq 2k-1.$$ 
 For $1\leq m\leq n\leq 2k$, we define in $\F_p[T]$ the continuants:$K_{m,n}=<v_mT,v_{m+1}T,\dots,v_{n}T>$ and $K_{1,0}=1$.
 \par For an integer $l\geq 1$, we introduce the subset $E(p,k,l)$ of hyperquadratic power series defined, via their continued fraction
expansion $\alpha=[a_1,a_2,\dots,a_l,\alpha_{l+1}]$, in the following way
$$(a_1,\dots,a_l)=(\lambda_1T,\dots,\lambda_lT) \quad \text{and }\quad \alpha^p=u_1K_{1,2k}\alpha_{l+1}+u_2K_{1,2k-1}$$
where $(\lambda_1,\lambda_2,\dots,\lambda_l,u_1,u_2)$ is an arbitrary $(l+2)$-tuple in $(\F_p^*)^{l+2}$.
These continued fractions have been studied in [4] and [5]. They all satisfy an algebraic equation of degree $p+1$ and therefore 
have an algebraic degree $d\leq p+1$. It has been proved that under a particular condition on this $(l+2)$-tuple,
the continued fraction has a very regular pattern which can be described. This condition is the following :
$$(*)\qquad [\lambda_1,\dots,\lambda_{l-1},\lambda_l-2k(u_1/u_2)(v_{2k}/v_1)]=k2^{-2k+1}\binom{2k}{k}u_2.$$
When condition $(*)$ is fullfilled, we say that the continued fraction expansion is perfect. To describe this expansion, 
we need to introduce the following sequence. Let $(A_m)_{m\geq 0}$ be the sequence in $\mathbb{F}_p[T]$ defined
recursively by (here the brackets denote the integer part of the rational) 
$$A_0=T \quad \text{ and}\quad
A_{m+1}=[A_m^{p}/(T^2-1)^k]\quad \text{ for}\quad m\geq 0.$$ 
Then for a perfect expansion $\alpha=[a_1,\dots,a_n,\dots] \in E(p,k,l)$, there exist a sequence $(\lambda_n)_{n\geq 1}$
in $\mathbb{F}_{p}^*$ and a sequence of non-negative integers $(i(n))_{n\geq 1}$ such that  
$$a_n=\lambda_nA_{i(n)}\quad \text{ for}\quad n\geq 1.$$
Both sequences $(\lambda_n)_{n\geq 1}$ and $(i(n))_{n\geq 1}$ are described in $[5]$. The sequence $(\lambda_n)_{n\geq 1}$ is defined
 recursively in a very sophisticated way involving the $(l+2)$-tuple appearing in the definition. The sequence $(i(n))_{n\geq 1}$ is 
very regular and it depends only upon the integers $k$ and $l$. Knowing this last sequence, the growth of the degree of the partial
 quotients is controlled and this allows to compute the quantity introduced at the end of Section 2. We have $t(\alpha)=(p-2k-1)/l$,
 for a perfect exansion $\alpha$ in $E(p,k,l)$. When condition $(*)$ is not fullfilled, the continued fraction has not been 
described. However, in such cases, it appears that the growth of the partial quotients is more important. Consequently this prevents
 the algebraic degree of such an element from being too small.
\par Now we turn to the root of $(RQE)$ (the reader may consult [6, p. 30-33]). We assume that $p$ is a prime such that $p\equiv 1 \mod 3$. Hence we 
can set $j=(p-1)/6$. We set $k=2j$ and we consider the $4j$-tuple $(v_1,\dots,v_{4j})$ introduced above. 
For $(\epsilon,\epsilon')\in (\F_p^*)^2$, we define in $\F_p[T][X]$ the following polynomial :
$$H(X)=K_{j+2,4j}X^{p+1}-\epsilon K_{j+1,4j}X^p+\epsilon'(K_{1,j}X+\epsilon K_{1,j-1}).$$
This polynomial has a unique irrational root $\alpha \in \F(p)$, with $\vert \alpha \vert \geq \vert T \vert$. This root 
belongs to $E(p,2j,3j)$ and the $(3j+2)$-tuple appearing in its definition is such that
$$u_2=(-1)^j\epsilon' ,\quad u_1=u_2\epsilon^{(-1)^{j+1}} \quad \text{and}\quad \lambda_i=v_{j+i}\epsilon^{(-1)^{i+1}}
\quad \text{for}\quad 1\leq i \leq 3j.$$
We want to prove that, for a particular choice of the pair $(\epsilon,\epsilon')\in (\F_p^*)^2$, this root of $H$
 is also the root of $(RQE)$. This will be obtained by proving that $P$ divides $H$, since $P$ has also only one root in $\F(p)$. 
Since the integer part of the root of $P$ is $(32/9)T$, while the integer part of the root of $H$ is $\epsilon v_{j+1}$, we need
 to have $\epsilon=32/(9v_{j+1})$. To guess the right value for $\epsilon'$, since the algebraic degree of the root of $P$ is 4 
and, for this root, the beginning of the continued fraction, observed by computer, matches to the pattern of a perfect expansion in 
$E(p,2j,3j)$, we will assume that condition $(*)$ holds. Hence we obtain
$$\epsilon'=\frac{(-16)^{j+1}}{3v_{j+1}\dbinom{4j}{2j}}[v_{j+1},v_{j+2},\dots,v_{4j-1},3v_{4j}/5].$$
Choosing the pair $(\epsilon,\epsilon')$ as indicated, the solution of $H$, and consequently the one of $P$ if $P$ divides $H$, 
will be defined by a perfect expansion in $E(p,2j,3j)$. Besides, this implies that we have $\deg(a_n)=(p^{i(n)}+2)/3$ and $t(\alpha)=2/3$. 
Moreover, in this particular case, the sequence $(i(n))_{n\geq 1}$ can be simply described. For $(m,n)\in\N^2$, we define
$${\bf{v}}_m[n]=\max \lbrace k\in \N \quad \text{such that} \quad m^k\vert n \rbrace$$
then we have $$i(n)={\bf{v}}_{4j+1}[4n-1]\quad \text{for} \quad n\geq 1.$$
\par At last, with these two values for the pair $(\epsilon,\epsilon')\in (\F_p^*)^2$, we need to prove that $P$ divides $H$. 
This can be checked directly by computer for a given prime $p$ with $p\equiv 1 \mod 3$ and we have done so for all primes up to 
$p=199$. Of course we conjecture that this is true for all primes $p$ with $p\equiv 1 \mod 3$.

\vskip 1 cm
\noindent Alain LASJAUNIAS
\newline Institut de Math\'ematiques de Bordeaux-CNRS-UMR 5251
\newline Talence 33405, France
\newline E-mail: Alain.Lasjaunias@math.u-bordeaux1.fr

\end{document}